\documentclass[11pt,a4paper]{amsart}
\usepackage{amsthm,amsmath,amsfonts,graphicx}
\usepackage[sort&compress,numbers]{natbib}
\usepackage[lmargin=37mm,rmargin=37mm,tmargin=35mm,bmargin=35mm]{geometry}

\newcommand{\spacing}[1]{
\renewcommand{\baselinestretch}{#1}
\setlength{\footnotesep}{\baselinestretch\footnotesep}}
\spacing{1.15}

\theoremstyle{plain}
\newtheorem{thm}{Theorem}[section]
\newtheorem{lem}[thm]{Lemma}
\newtheorem{cor}[thm]{Corollary}
\newtheorem{claim}[thm]{Claim}

\newcommand{\COR}[2]{\begin{cor}\label{cor:#1}#2\end{cor}}
\newcommand{\THM}[2]{\begin{thm}\label{thm:#1}#2\end{thm}}
\newcommand{\LEM}[2]{\begin{lem}\label{lem:#1}#2\end{lem}}
\newcommand{\LEMMA}[3]{\begin{lem}[#1]\label{lem:#2}#3\end{lem}}
\newcommand{\SECT}[2]{\section{#2}\label{sec:#1}}
\newcommand{\EQN}[2]{\begin{equation}\label{eqn:#1}#2\end{equation}}
\newcommand{\PROOF}[1]{\begin{proof}#1\end{proof}}
\newcommand{\BLAH}[2]{\begin{proof}[#1]#2\end{proof}}

\newcommand{\corref}[1]{Corollary~\ref{cor:#1}}
\newcommand{\thmref}[1]{Theorem~\ref{thm:#1}}
\newcommand{\lemref}[1]{Lemma~\ref{lem:#1}}
\newcommand{\secref}[1]{Section~\ref{sec:#1}}
\newcommand{\twosecref}[2]{Sections~\ref{sec:#1} and \ref{sec:#2}}

\newcommand{\eqnref}[1]{Equation~\eqref{eqn:#1}}
\newcommand{\twoeqnref}[2]{Equations~\eqref{eqn:#1} and \eqref{eqn:#2}}
\newcommand{\figref}[1]{Figure~\ref{fig:#1}}

\newcommand{\FLOOR}[1]{\ensuremath{\protect\left\lfloor#1\right\rfloor}}
\newcommand{\floor}[1]{\ensuremath{\protect\lfloor#1\rfloor}}
\newcommand{\ceil}[1]{\ensuremath{\protect\lceil#1\rceil}}
\newcommand{\twothirds}{\ensuremath{\protect\tfrac{2}{3}}}
\newcommand{\CP}[2]{\ensuremath{#1\,\square\,#2}}
\newcommand{\CCP}[3]{\ensuremath{\protect#1{\,\square\,}#2{\,\square\,}\cdots{\,\square\,}#3}}

\begin{document}

\title[Metric Dimension of Cartesian Products]{On the Metric Dimension of\\ Cartesian Products of Graphs}

\thanks{Research supported by projects MCYT-FEDER-BFM2003-00368, Gen-Cat-2001SGR00224, MCYT-HU2002-0010, MTM-2004-07891-C02-01, MEC-SB2003-0270, MCYT-FEDER BFM2003-00368, and Gen.\ Cat 2001SGR00224.}

\author[C{\'a}ceres]{Jos{\'e} C{\'a}ceres}
\address{Departamento de Estad{\'\i}stica y Matem{\'a}tica Aplicada, Universidad de Almer{\'\i}a, Almer{\'\i}a, Spain}
\email{jcaceres@ual.es}

\author[Hernando]{Carmen Hernando}
\address{Departament de Matem{\`a}tica Aplicada I, Universitat Polit{\`e}cnica de Catalunya, Barcelona, Spain}
\email{carmen.hernando@upc.edu}

\author[Mora]{Merc{\`e} Mora}
\address{Departament de Matem{\`a}tica Aplicada II, Universitat Polit{\`e}cnica de Catalunya, Barcelona, Spain}
\email{merce.mora@upc.edu}

\author[Pelayo]{Ignacio M. Pelayo}
\address{Departament de Matem{\`a}tica Aplicada III, Universitat Polit{\`e}cnica de Catalunya, Barcelona, Spain}
\email{ignacio.m.pelayo@upc.edu}

\author[Puertas]{Mar{\'\i}a L.~Puertas}
\address{Departamento de Estad{\'\i}stica y Matem{\'a}tica Aplicada, Universidad de Almer{\'\i}a, Almer{\'\i}a, Spain}
\email{mpuertas@ual.es}

\author[Seara]{Carlos Seara}
\address{Departament de Matem{\`a}tica Aplicada II, Universitat Polit{\`e}cnica de Catalunya, Barcelona, Spain}
\email{carlos.seara@upc.edu}

\author[Wood]{David R.~Wood}
\address{Departament de Matem{\`a}tica Aplicada II, Universitat Polit{\`e}cnica de Catalunya, Barcelona, Spain}
\email{david.wood@upc.edu}



\begin{abstract}
A set of vertices $S$ \emph{resolves} a graph $G$ if every vertex
is uniquely determined by its vector of distances to the vertices
in $S$. The \emph{metric dimension} of $G$ is the minimum
cardinality of a resolving set of $G$. This paper studies the
metric dimension of cartesian products $G\,\square\,H$. We prove that
the metric dimension of $G\,\square\,G$ is tied in a strong sense to the
minimum order of a so-called doubly resolving set in $G$. Using
bounds on the order of doubly resolving sets, we establish bounds
on $G\,\square\,H$ for many examples of $G$ and $H$. One of our main
results is a family of graphs $G$ with bounded metric dimension
for which the metric dimension of $G\,\square\,G$ is unbounded.
\end{abstract}

\maketitle

\SECT{Introduction}{Introduction}

A set of vertices $S$ \emph{resolves} a graph $G$ if every vertex
of $G$ is uniquely determined by its vector of distances to the
vertices in $S$. This paper undertakes a general study of
resolving sets in cartesian products of graphs.

All the graphs considered are finite, undirected, simple, and
connected\footnote{The results can easily be generalised to
disconnected graphs; we omit the details.}. The vertex set and
edge set of a graph $G$ are denoted by $V(G)$ and $E(G)$. The
distance between vertices $v,w\in V(G)$ is denoted by $d_G(v,w)$,
or $d(v,w)$ if the graph $G$ is clear from the context. A vertex
$x\in V(G)$ \emph{resolves} a pair of vertices $v,w\in V(G)$ if
$d(v,x)\ne d(w,x)$. A set of vertices $S\subseteq V(G)$
\emph{resolves} $G$, and $S$ is a \emph{resolving set} of $G$, if
every pair of distinct vertices of $G$ are resolved by some vertex
in $S$. A resolving set $S$ of $G$ with the minimum cardinality is
a \emph{metric basis} of $G$, and $|S|$ is the \emph{metric
dimension} of $G$, denoted by $\beta(G)$.

The \emph{cartesian product} of graphs $G$ and $H$, denoted by
\CP{G}{H}, is the graph with vertex set $V(G)\times
V(H):=\{(a,v):a\in V(G),v\in V(H)\}$, where $(a,v)$ is adjacent to
$(b,w)$ whenever $a=b$ and $\{v,w\}\in E(H)$, or $v=w$ and
$\{a,b\}\in E(G)$. Where there is no confusion the vertex $(a,v)$
of \CP{G}{H} will be written $av$. Observe that if $G$ and $H$ are
connected, then \CP{G}{H} is connected. In particular, for all
vertices $av,bw$ of \CP{G}{H} we have
$d(av,bw)=d_G(a,b)+d_H(v,w)$. Assuming isomorphic graphs are
equal, the cartesian product is associative, and
$\CCP{G_1}{G_2}{G_d}$ is well-defined.

Resolving sets in general graphs were first defined by
\citet{HM-AC76} and \citet{Slater75}, although as we shall see,
resolving sets in hypercubes were studied earlier under the guise
of a coin weighing problem \citep{Lindstrom64, ER63, CM-CJM66,
Lindstrom66, AKV-FOCS96, Pippenger77, SS-AMM63, FS-AMM67,
Cantor-CJM64, KLT00, Lindstrom-CMB65, GN-AMM95, Mow77}. Resolving
sets have since been widely investigated \citep{SSH02, SZ-CN04,
SZ-IJMMS04, CZ-CN03, ST-MOR04, BCDZ-MB03, PZ02, CO01, SS01,
CEJO-DAM00, CPZ-CMA00, Sooryanarayana98, KRR-DAM96, Yushmanov87,
SZ-CMJ04, CGH, Slater-JMPS88, SZ-CMJ03}, and arise in many diverse
areas including network discovery and verification
\citep{BEEHHMS-WG05}, robot navigation \citep{KRR-DAM96, SSH02},
connected joins in graphs \citep{ST-MOR04}, and strategies for the
Mastermind game \citep{BG, Greenwell-JRM99, Goddard-JCMCC03,
Goddard-JCMCC04, KLT00, Chvatal-Comb83}.


Part of our motivation for studying the metric dimension of
cartesian products is that in two of the above-mentioned
applications, namely Mastermind and coin weighing, the graphs that
arise are in fact cartesian products. These connections are
explained in \twosecref{Coin}{Hamming} respectively.

The main contributions of this paper are based on the notion of
doubly resolving sets, which are introduced in
\secref{DoublyResolving}. We prove that the minimum order of a
doubly resolving set in a graph $G$ is tied in a strong sense to
$\beta(\CP{G}{G})$. Thus doubly resolving sets are essential in
the study of metric dimension of cartesian products. We then give
a number of examples of bounds on the metric dimension of
cartesian products through doubly resolving sets. In particular,
Sections~\ref{sec:CompleteGraphs}, \ref{sec:Hamming},
\ref{sec:Paths}, \ref{sec:Cycles}, and \ref{sec:Trees}
respectively study complete graphs, Hamming graphs, paths and
grids, cycles, and trees. One of our main results here is a family
of (highly connected) graphs with bounded metric dimension for
which the metric dimension of the cartesian product is unbounded.

\SECT{Coin}{Coin Weighing and Hypercubes}

The \emph{hypercube} $Q_n$ is the graph whose vertices are the
$n$-dimensional binary vectors, where two vertices are adjacent if
they differ in exactly one coordinate. It is well known that
\begin{equation*}
Q_n=\underbrace{\CCP{K_2}{K_2}{K_2}}_n.
\end{equation*}
It is easily seen that $\beta(Q_n)\leq n$ (see \eqnref{Grid}). The
first case when this bound is not tight is $n=5$. A laborious
calculation verifies that $Q_5$ is resolved by the $4$-vertex set
$\{00000,00011,00101,01001\}$. We have determined $\beta(Q_n)$ for
small values of $n$ by computer search.

\medskip
\begin{center}
\begin{tabular}{c|ccccccccc} \hline
 $n$      &  $2$     & $3$  &  $4$    &  $5$    & $6$  &  $7$   & $8$ & $10$& $15$\\  \hline
 $\beta(Q_n)$   & $2 $ & $3$  &  $4$    &  $4$ &    $5$  & $6$  &   $6$ & $\leq7$& $\leq10$\\ \hline
\end{tabular}
\end{center}
\medskip

The asymptotic value of $\beta(Q_n)$ turns out to be related to
the following coin weighing problem first posed by
\citet{SS-AMM63}. (See \citep{GN-AMM95} for a survey on various
coin weighing problems.)\ Given $n$ coins, each with one of two
distinct weights, determine the weight of each coin with the
minimum number of weighings. We are interested in the static
variant of this problem, where the choice of sets of coins to be
weighed is determined in advance. Weighing a set $S$ of coins
determines how many light (and heavy) coins are in $S$, and no
further information. It follows that the minimum number of
weighings differs from $\beta(Q_n)$ by at most one
\citep{ST-MOR04, KLT00}. A lower bound on the number of weighings
by \citet{ER63} and an upper bound by \citet{Lindstrom64} imply
that
\begin{equation*}
\lim_{n\rightarrow\infty}\,\beta(Q_n)\cdot\frac{\log n}{n}=2,
\end{equation*}
where, as always in this paper, logarithms are binary. Note that
Lindstrom's proof is constructive. He gives an explicit scheme of $2^k-1$
weighings that suffice for $k\cdot2^{k-1}$ coins. 

\SECT{Projections}{Projections}

Let $S$ be a set of vertices in the cartesian product \CP{G}{H} of
graphs $G$ and $H$. The \emph{projection} of $S$ onto $G$ is the
set of vertices $a\in V(G)$ for which there exists a vertex $av\in
S$. Similarly, the \emph{projection} of $S$ onto $H$ is the set of
vertices $v\in V(H)$ for which there exists a vertex $av\in S$. A
\emph{column} of \CP{G}{H} is the set of vertices $\{av:v\in
V(H)\}$ for some vertex $a\in V(G)$, and a \emph{row} of \CP{G}{H}
is the set of vertices $\{av:a\in V(G)\}$ for some vertex $v\in
V(H)$. Observe that each row induces a copy of $G$, and each
column induces a copy of $H$. This terminology is consistent with
a representation of \CP{G}{H} by the points of the
$|V(G)|\times|V(H)|$ grid.

\LEM{Projections}{Let $S\subseteq V(\CP{G}{H})$ for graphs $G$ and
$H$. Then every pair of vertices in a fixed row of \CP{G}{H} are
resolved by $S$ if and only if the projection of $S$ onto $G$
resolves $G$. Similarly, every pair of vertices in a fixed column
of \CP{G}{H} are resolved by $S$ if and only if the projection of
$S$ onto $H$ resolves $H$.}

\PROOF{Consider two vertices $av$ and $aw$ in a common column. For
every other vertex $bx$ of \CP{G}{H}, we have
$d(av,bx)-d(aw,bx)=d_H(v,x)-d_H(w,x)$. Thus $d(av,bx)\ne d(aw,bx)$
if and only if $d_H(v,x)\ne d_H(w,x)$. That is, $av$ and $aw$ are
resolved by $bx$ if and only if $v$ and $w$ are resolved by $x$ in
$H$. Hence $av$ and $aw$ are resolved by $S$ if and only if $v$
and $w$ are resolved by the projection of $S$ onto $H$. We have
the analogous result for the projection onto $G$ by symmetry.}

\COR{Projections}{For all graphs $G$ and $H$, and for every
resolving set $S$ of \CP{G}{H}, the projection of $S$ onto $G$
resolves $G$, and the projection of $S$ onto $H$ resolves $H$. In
particular, $\beta(\CP{G}{H})\geq\max\{\beta(G),\beta(H)\}$.\qed}

\SECT{DoublyResolving}{Doubly Resolving Sets}

Many of the results that follow are based on the following
definitions. Let $G\ne K_1$ be a graph. Two vertices $v,w\in V(G)$
are \emph{doubly resolved} by $x,y\in V(G)$ if
\begin{equation*}
d(v,x)-d(w,x)\ne d(v,y)-d(w,y).
\end{equation*}

A set of vertices $S\subseteq V(G)$ \emph{doubly resolves} $G$,
and $S$ is a \emph{doubly resolving set}, if every pair of
distinct vertices $v,w\in V(G)$ are doubly resolved by two
vertices in $S$. Every graph with at least two vertices has a
doubly resolving set. Let $\psi(G)$ denote the minimum cardinality
of a doubly resolving set of a graph $G\ne K_1$. Note that if
$x,y$ doubly resolves $v,w$ then $d(v,x)-d(w,x)\ne0$ or
$d(v,y)-d(w,y)\ne0$, and at least one of $x$ and $y$ (singly)
resolves $v,w$. Thus a doubly resolving set is also a resolving
set, and
\begin{equation*}
\beta(G)\leq\psi(G).
\end{equation*}
Our interest in doubly resolving sets is based on the following upper bound.

\THM{DoublyResolving}{For all graphs $G$ and $H\ne K_1$,
\begin{equation*}
\beta(\CP{G}{H})\leq\beta(G)+\psi(H)-1.
\end{equation*}}

\PROOF{Let $S$ be a metric basis of $G$.
Let $T$ be a doubly resolving set of $H$ with $|T|=\psi(H)$.
Fix vertices $s\in S$ and $t\in T$.
Let \begin{equation*}
X:=\{sv:v\in T\}\cup\{at:a\in S\}.
\end{equation*}
Observe that $|X|=|S|+|T|-1$. To prove that $X$ resolves
\CP{G}{H}, consider two vertices $av$ and $bw$ of \CP{G}{H}. By
\lemref{Projections}, if $a=b$ then $av$ and $bw$ are resolved
since the projection of $X$ onto $H$ is $T$. Similarly, if $v=w$
then $av$ and $bw$ are resolved since the projection of $X$ onto
$G$ is $S$. Now assume that $a\ne b$ and $v\ne w$. Since $T$ is
doubly resolving for $H$, there are two vertices $x,y\in T$ such
that
\begin{equation*}
d_H(v,x)-d_H(w,x)\ne d_H(v,y)-d_H(w,y).
\end{equation*}
Thus for at least one of $x$ and $y$, say $x$, we have
\begin{equation*}
d_H(v,x)-d_H(w,x)\ne d_G(b,s)-d_G(a,s).
\end{equation*}
Hence \begin{equation*}
d(av,sx)=d_G(a,s)+d_H(v,x)\ne d_G(b,s)+d_H(w,x)=d(bw,sx).
\end{equation*}
That is, $sx\in X$ resolves $av$ and $bw$.}

The relationship between resolving sets of cartesian products and
doubly resolving sets is strengthened by the following lower
bound.

\LEM{TwoProjections}{Suppose that $S$ resolves \CP{G}{G} for some graph $G$.
Let $A$ and $B$ be the two projections of $S$ onto $G$.
Then $A\cup B$ doubly resolves $G$. In particular,
\begin{equation*}
\beta(\CP{G}{G})\geq\tfrac{1}{2}\psi(G).
\end{equation*}}

\PROOF{For any two vertices $v,w\in V(G)$, there is a vertex
$pq\in S$ that resolves $vw,wv$. That is, $d(vw,pq)\ne d(wv,pq)$.
Thus $d(v,p)+d(w,q)\ne d(w,p)+d(v,q)$, which implies
$d(v,p)-d(w,p)\ne d(v,q)-d(w,q)$. Thus $p,q$ doubly resolves $v,w$
in $G$. Now $p\in A$ and $q\in B$. Hence $A\cup B$ doubly resolves
$G$. If, in addition, $S$ is a metric basis of \CP{G}{G}, then
$\psi(G)\leq|A\cup B|\leq|A|+|B|\leq2|S|=2\cdot\beta(\CP{G}{G})$.}

Observe that \thmref{DoublyResolving} and \lemref{TwoProjections} prove
that $\beta(\CP{G}{G})$ is always within a constant factor of $\psi(G)$. In particular,
\begin{equation}
\label{eqn:Summarise}
\tfrac{1}{2}\psi(G)
\leq\beta(\CP{G}{G})
\leq\psi(G)+\beta(G)-1
\leq2\psi(G)-1.
\end{equation}
Thus doubly resolving sets are essential in the study of the
metric dimension of cartesian products.

A natural candidate for a resolving set of \CP{G}{G} is $S\times
S$ for a well chosen set $S\subseteq V(G)$. It follows from
\lemref{TwoProjections} and the proof technique employed in
\thmref{DoublyResolving} that $S\times S$ resolves \CP{G}{G} if
and only if $S$ doubly resolves $G$.

Now consider the following elementary bound on $\psi(G)$.

\LEM{PsiWeakBound}{For every graph $G$ with $n\geq3$ vertices we have $\psi(G)\leq n-1$.}

\PROOF{Clearly $G$ has a vertex $x$ of degree at least two. Let
$S:=V(G)\setminus\{x\}$. To prove that $S$ doubly resolves $G$,
consider two vertices $u,v\in V(G)$. If both $u,v\in S$, then the
pair $u,v$ doubly resolves itself. Otherwise, without loss of
generality, $u\in S$ and $v=x$. Since $\deg(x)\geq2$, there is a
neighbour $y\ne u$ of $x$. Now $d(u,u)-d(v,u)\leq 0-1=-1$ and
$d(u,y)-d(v,y)\geq 1-1=0$. Thus $u,y\in S$ doubly resolve $u,v$.
Hence $S$ doubly resolves $G$.}

Note that if $G$ is a graph with $n\geq3$ vertices, then
\thmref{DoublyResolving} and \lemref{PsiWeakBound} imply that
$\beta(\CP{G}{H})\leq\beta(H)+n-2$ for every graph $H$.

\SECT{CompleteGraphs}{Complete Graphs}

Let $K_n$ denote the complete graph on $n\geq1$ vertices. It is
well known \citep{KRR-DAM96, CEJO-DAM00} that for every $n$-vertex
graph $G$,
\EQN{BetaCompleteGraph}{\beta(G)=n-1\,\Longleftrightarrow\,G=K_n.}

\LEM{PsiCompleteGraph}{For all $n\geq2$ we have $\psi(K_n)=\max\{n-1,2\}$.}

\PROOF{Since $\psi(G)\geq2$ for every graph $G\ne K_1$, we have $\psi(K_2)=2$.
Now suppose that $n\geq3$. By \lemref{PsiWeakBound}, $\psi(K_n)\leq n-1$.
Conversely, $\psi(K_n)\geq\beta(K_n)=n-1$ by \eqnref{BetaCompleteGraph}.}


\thmref{DoublyResolving} and \lemref{PsiCompleteGraph} imply that every graph $G$ satisfies
\EQN{CompleteTimesGraph}{\beta(\CP{K_n}{G})\leq\beta(G)+\max\{n-2,1\}.}
In certain cases, this result can be improved as follows.

\LEM{CompleteGraphProduct}{For every graph $G$ and for all $n\geq1$,
\begin{equation*}
\beta(\CP{K_n}{G})\leq\max\{n-1,2\cdot\beta(G)\}.
\end{equation*}}

\PROOF{Let $S$ be a metric basis of $G$. Fix a vertex $r$ of
$K_n$. As illustrated in \figref{CompleteTimesGraph}, there is a
set $T$ of $\max\{n-1,2|S|\}$ vertices of \CP{K_n}{G} such that:
\begin{enumerate}
\item[\textup{(a)}] for all vertices $a\in V(K_n)\setminus\{r\}$, there is
at least one vertex $x\in S$ for which $ax\in T$, and
\item[\textup{(b)}] for all $x\in S$, there are at least two vertices
$a,b\in V(K_n)$ for which $ax\in T$ and $bx\in T$.
\end{enumerate}

\begin{figure}[!ht]
\begin{center}
\includegraphics{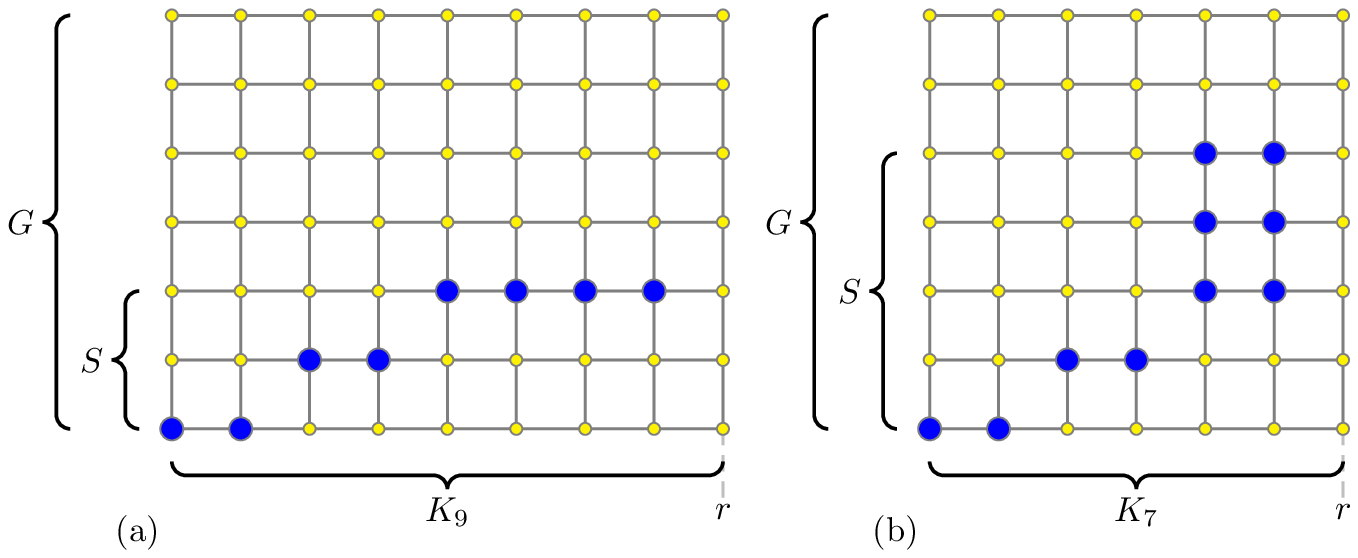}
\end{center}
\caption{\label{fig:CompleteTimesGraph} The resolving set $T$ of \CP{K_n}{G} in \lemref{CompleteGraphProduct}: (a) $n-1\geq2\beta(G)$ and (b) $n-1\leq2\beta(G)$.}
\end{figure}

To prove that $T$ resolves \CP{K_n}{G}, consider two vertices $av$ and $bw$ of \CP{K_n}{G}. If $v=w$, then since the projection of $T$ onto $G$ is the resolving set $S$, by \lemref{Projections}, $av$ and $bw$ are resolved by $T$. Now suppose that $v\ne w$. Then there is a vertex $x\in S$ that resolves  $v$ and $w$ in $G$. Hence $d_G(v,x)<d_G(w,x)$ without loss of generality. By (b) there are distinct vertices $c,d\in V(K_n)$ for which $cx\in T$ and $dx\in T$. If $c\ne a$ and $c\ne b$, then
\begin{equation*}
d(av,cx)=d_G(v,x)+1<d_G(w,x)+1=d(bw,cx));
\end{equation*}
that is, $cx$ resolves $av$ and $bw$ in \CP{K_n}{G}. Similarly, if  $d\ne a$ and $d\ne b$, then $dx$ resolves $av$ and $bw$. Otherwise $c=a$ or $c=b$, and $d=a$ or $d=b$. Since $c\ne d$, without loss of generality $c=a$ and $d=b$. Then \begin{equation*}
d(av,cx)=d_G(v,x)<d_G(w,x)<d_G(w,x)+1=d(bw,cx),
\end{equation*}
and again $cx$ resolves $av$ and $bw$ in \CP{K_n}{G}.}

When is $n$ is large in comparison with $\beta(G)$ we know $\beta(\CP{K_n}{G})$ exactly.

\THM{BigClique}{For every graph $G$ and for all $n\geq2\cdot\beta(G)+1$,
\begin{equation*}
\beta(\CP{K_n}{G})=n-1.
\end{equation*}}

\PROOF{The lower bound $\beta(\CP{K_n}{G})\geq n-1$ follows from \corref{Projections} and \eqnref{BetaCompleteGraph}. The upper bound $\beta(\CP{K_n}{G})\leq n-1$ is a special case of \lemref{CompleteGraphProduct}.}

\SECT{Hamming}{Mastermind and Hamming Graphs}

\emph{Mastermind} is a game for two players, the \emph{code setter} and the \emph{code breaker}\footnote{\citet{Chvatal-Comb83} referred to the code setter and code breaker as S.F.\ and P.G.O.M.\ (in honour of P.E.).}. The code setter chooses a secret vector $s=[s_1,s_2,\dots,s_n]\in\{1,2,\dots,k\}^n$. The task of the code breaker is to infer the secret vector by a series of questions, each a vector $t=[t_1,t_2,\dots,t_n]\in\{1,2,\dots,k\}^n$. The code setter answers with two integers, first being the number of positions in which the secret vector and the question agree, denoted by $a(s,t)=|\{i:s_i=t_i,1\leq i\leq n\}|$. The second integer $b(s,t)$ is the maximum of $a(\tilde{s},t)$, where $\tilde{s}$ ranges over all permutations of $s$.

In the commercial version of the game, $n=4$ and $k=6$. The secret vector and each question is represented by four pegs each coloured with one of six colours. Each answer is represented by $a(s,t)$ black pegs, and $b(s,t)-a(s,t)$ white pegs. \citet{Knuth-JRM76} showed that four questions suffice to determine $s$ in this case. Here the code breaker may determine each question in response to the previous answers. \emph{Static mastermind} is the variation in which all the questions must be supplied at once. Let $g(n,k)$ denote the maximum, taken over all vectors $s$, of the minimum number of questions required to determine $s$ in this static setting.

The \emph{Hamming graph} $H_{n,k}$ is the cartesian product of cliques \begin{equation*}
H_{n,k}=\underbrace{\CCP{K_k}{K_k}{K_k}}_n.
\end{equation*}
Note that the hypercube $Q_n=H_{n,2}$. The vertices of $H_{n,k}$ can be thought of as vectors in $\{1,2,\dots,k\}^n$, with two vertices being adjacent if they differ in precisely one coordinate. Thus the distance $d_H(v,w)$ between two vertices $v$ and $w$ is the number of coordinates in which their vectors differ. That is,
\begin{equation*}
d_H(v,w)=n-a(v,w).
\end{equation*}

Suppose for the time being that we remove the second integer $b(s,t)$ from the answers given by the code setter in the static mastermind game. Let $f(n,k)$ denote the maximum, taken over all vectors $s$, of the minimum number of questions required to determine $s$ without $b(s,t)$ in the answers. For the code breaker to correctly infer the secret vector $s$ from a set of questions $T$, $s$ must be uniquely determined by the values $\{a(s,t):t\in T\}$. Equivalently, for any two vertices $v$ and $w$ of $H_{n,k}$, there is a $t\in T$ for which $a(v,t)\ne a(w,t)$; that is, the distances $d_H(v,t)\ne d_H(w,t)$. Hence the secret vector can be inferred if and only if $T$ resolves $H_{n,k}$. Thus
\begin{equation*}
g(n,k)\leq f(n,k)=\beta(H_{n,k}).
\end{equation*}
\citet{Chvatal-Comb83} proved the upper bound
\begin{equation*}
\beta(H_{n,k})=f(n,k)\leq(2+\epsilon)n\,\frac{1+2\log k}{\log n-\log k}
\end{equation*}
for large $n>n(\epsilon)$ and small $k<n^{1-\epsilon}$. For $k\in\{3,4\}$, improvements to the constant in the above upper bound are stated without proof by \citet{KLT00}.
%
%
They also state that a `straightforward generalisation' of the lower bound on $\beta(Q_n)$ by \citet{ER63} gives for large $n$,
\begin{equation*}
\beta(H_{n,k})\geq g(n,k)\geq(2+o(1))\frac{n\log k}{\log n}.
\end{equation*}

Here we study $\beta(H_{n,k})$ for large values of $k$ rather that for large values of $n$. A similar approach is take by \citet{Goddard-JCMCC03, Goddard-JCMCC04} for static Mastermind, who proved that $g(2,k)=\ceil{\frac{2}{3}k}$ and $g(3,k)=k-1$. Our contribution is to determine the exact value of $\beta(H_{2,k})$. We show that for all $k\geq1$,
\EQN{TwoDimHamming}{\beta(H_{2,k})=\FLOOR{\tfrac{2}{3}(2k-1)}.}
\eqnref{TwoDimHamming} is a special case (with $m=n=k$) of the following more general result.

\THM{TwoCliquesProduct}{For all $n\geq m\geq1$ we have
\begin{equation*}
\beta(\CP{K_n}{K_m})=
\begin{cases}
    \FLOOR{\tfrac{2}{3}(n+m-1)} & \textup{, if }m\leq n\leq 2m-1\\
    n-1                         & \textup{, if }n\geq 2m-1.
\end{cases}
\end{equation*}}

\begin{figure}
\begin{center}
\includegraphics{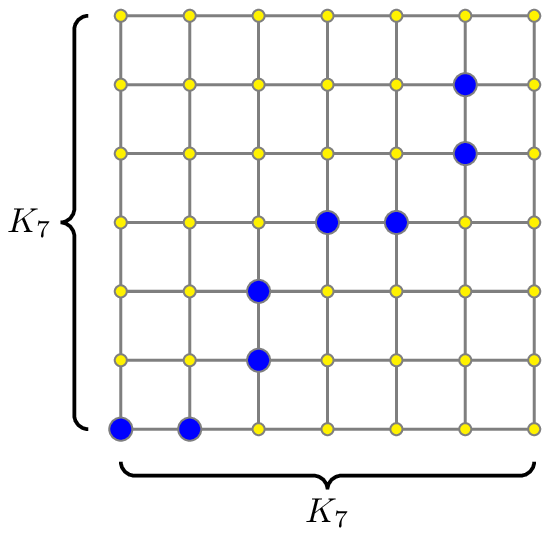}
\end{center}
\caption{\label{fig:CompleteComplete} Resolving set of \CP{K_7}{K_7} with one empty row, one empty column, and no lonely vertex.}
\end{figure}

Note that two vertices of \CP{K_n}{K_m} are adjacent if and only if they are in a common row or column. Otherwise they are at distance two. Fix a set $S$ of vertices of \CP{K_n}{K_m}. With respect to $S$, a row or column is \emph{empty} if it contains no vertex in $S$, and a vertex $v\in S$ is \emph{lonely} if $v$ is the only vertex of $S$ in its row and in its column. As illustrated in \figref{CompleteComplete}, we have the following characterisation of resolving sets in  \CP{K_n}{K_m}.

\LEM{TwoCliquesCharac}{For $m,n\ge 2$, a set $S$ of vertices resolves \CP{K_n}{K_m} if and only if:
\begin{enumerate}
\item[\textup{(a)}] there is at most one empty row and at most one empty column,
\item[\textup{(b)}] there is at most one lonely vertex, and
\item[\textup{(c)}] if there is an empty row and an empty column, then there is no lonely vertex.
\end{enumerate}}

\PROOF{($\Longrightarrow$) First suppose that $S$ resolves \CP{K_n}{K_m}. By \corref{Projections}, the projections of $S$ respectively resolve $K_m$ and $K_n$. By \eqnref{BetaCompleteGraph}, there is at most one empty row and at most one empty column. Thus (a)  holds.

Suppose on the contrary that $v$ and $w$ are two lonely vertices in $S$. Thus $v$ and $w$ are in distinct rows and  columns, and no other vertex of $S$ is in a row or column that contains $v$ or $w$. Let $x$ be the vertex in the row of $v$ and the column of $w$. Let $y$ be the vertex in the column of $v$ and the row of $w$. Then $d(x,v)=d(y,v)=1$, $d(x,w)=d(y,w)=1$, and $d(x,u)=d(y,u)=2$ for every vertex $u\in S\setminus\{v,w\}$. Thus $S$ does not resolve $x$ and $y$. This contradiction proves that $S$ satisfies (b).

Finally, suppose that there is an empty row, an empty column, and a lonely  vertex $v\in S$. Let $x$ be the vertex in the row of $v$ and in the empty column. Let $y$ be the vertex in the column of $v$ and in the empty row. We have
$d(x,v)=d(y,v)=1$, and $d(x,u)=d(y,u)=2$ for every vertex $u\in S\setminus\{v\}$. Thus $S$ does not resolve $x$ and $y$. This contradiction proves that $S$ satisfies (c).

($\Longleftarrow$) Now suppose that $S$ is a set of vertices satisfying (a), (b) and (c). We will prove that $S$ resolves any two vertices $x$ and $y$. If $x\in S$, then $x$ resolves $x,y$. If $y\in S$, then $y$ resolves $x,y$. Now suppose that $x\not\in S$ and $y\not\in S$.

If $x$ and $y$ are in the same row, then at least one of the columns of $x$ and $y$ contains a vertex $v\in S$. Suppose $v$ is in the column of $x$. Thus $d(x,v)=1$ and $d(y,v)=2$, and $v$ resolves $x,y$.
Similarly, if $x$ and $y$ are in the same column, then some $v\in S$ resolves $x,y$.

Suppose now that $x$ and $y$ are in distinct rows and columns. Then there is a vertex of $S$ in the column of $x$ or in the column of $y$. Suppose $v\in S$ is in the column of $x$. If $v$ is not in the row of $y$, $d(x,v)=1\ne 2=d(y,v)$, and $v$ resolves $x,y$. If $v$ is in the row of $y$, by (b) and (c), at least one of the vertices in the rows and columns of $x$ and $y$, but not in the intersection of two of them, is in $S$. This vertex resolves $x$ and $y$.}


\LEM{Move}{For all $n,m\geq3$, if $S$ resolves \CP{K_n}{K_m}, then there exists a resolving set $S^*$ of \CP{K_n}{K_m} such that $|S^*|\leq|S|$, and $S$ contains two vertices $v$ and $w$ in the same row or column, such that $v$ and $w$ are the only vertices in $S^*$ in the row(s) and column(s) that contain $v$ and $w$.}

\PROOF{By \lemref{TwoCliquesCharac}, there are two vertices $v,w\in S$ in the same row or column. By symmetry, we can suppose that $v$ and $w$ are in the same row. If $v$ and $w$ are the only vertices in $S^*$ in the row and columns that contain $v$ and $w$, then we are done. Otherwise there is a vertex $x\in S$ in the row or columns that contain $v$ and $w$. It suffices to prove that $x$ can be deleted from $S$, or replaced in $S$ by some other vertex not in the row or columns that contain $v$ and $w$, such that $S$ still satisfies the conditions of \lemref{TwoCliquesCharac}, and thus resolves  \CP{K_n}{K_m}. We can then repeat this step to obtain the desired set $S^*$.

First suppose that $x$ is in the same row as $v$ and $w$. If all
the vertices of the column of $x$ are in $S$, then delete $x$ from
$S$; clearly $S$ still satisfies the conditions of
\lemref{TwoCliquesCharac}. Otherwise, let $y$ be a vertex not in
$S$ such that $y$ is in the column containing $x$, and if $x$ is
the only vertex in its column that is in $S$, then $y$ is in a row
that contains at least one vertex of $S$. This is always possible,
since $S$ satisfies condition (a). Then
$(S\setminus\{x\})\cup\{y\}$ satisfies the conditions of
\lemref{TwoCliquesCharac}.

Now suppose that $x$ is in the column of $v$ or $w$. If every vertex in the row containing $x$ is in $S$, then delete $x$ from $S$; clearly $S$ still satisfies the the conditions of \lemref{TwoCliquesCharac}. Otherwise, proceeding as in the preceding case, let $y$ be a vertex in the same row as $x$, but not in the columns of $v$ and $w$, such that there is at least one other vertex of $S$ in the row or column that contains $y$. Then $(S\setminus\{x\})\cup\{y\}$ satisfies the conditions of \lemref{TwoCliquesCharac}. This completes the proof.}


\LEM{TwoCliquesProduct}{For all $n,m\geq3$,
\begin{equation*}
\beta(\CP{K_n}{K_m})
=2+\min\{\beta(\CP{K_{n-2}}{K_{m-1}}),\beta(\CP{K_{n-1}}{K_{m-2}})\}.
\end{equation*}}

\PROOF{We first prove that
\EQN{TwoCliquesProductA}{\beta(\CP{K_n}{K_m}) \leq 2+\min\{\beta(\CP{K_{n-2}}{K_{m-1}}),\beta(\CP{K_{n-1}}{K_{m-2}})\}.}
Without loss of generality $\beta(\CP{K_{n-2}}{K_{m-1}})\leq\beta(\CP{K_{n-1}}{K_{m-2}})$. Let $S$ be a metric basis of \CP{K_{n-2}}{K_{m-1}}. Construct $S'\subseteq V(\CP{K_n}{K_m})$ from $S$ by adding two new vertices that are positioned in one new row and in two new columns. The number of empty rows, empty columns, and lonely vertices is the same in $S$ and $S'$. Since $S$ resolves \CP{K_{n-2}}{K_{m-1}}, $S'$ resolves \CP{K_n}{K_m} by \lemref{TwoCliquesCharac}.
Thus $\beta(\CP{K_n}{K_m})\leq|S'|=|S|+2=2+\beta(\CP{K_{n-2}}{K_{m-1}})$, which implies \eqref{eqn:TwoCliquesProductA}. It remains to prove that
\EQN{TwoCliquesProductB}{\min\{\beta(\CP{K_{n-2}}{K_{m-1}}),\beta(\CP{K_{n-1}}{K_{m-2}})\}\leq
\beta(\CP{K_n}{K_m})-2.}
Let $S$ be a metric basis of \CP{K_n}{K_m}. By \lemref{Move}, we can assume that $S$ contains two vertices $v$ and $w$ in the same row or column, such that $v$ and $w$ are the only vertices in $S$ in the row(s) and column(s) that contain $v$ and $w$. Without loss of generality, $v$ and $w$ are in the same row. Construct $S'\subseteq V(\CP{K_{n-2}}{K_{m-1}})$ from $S$ by deleting the row containing $v$ and $w$, and by deleting the two columns containing $v$ and $w$. The number of empty rows, empty columns, and lonely vertices is the same in $S$ and $S'$. Since $S$ resolves \CP{K_n}{K_m}, $S'$ resolves \CP{K_{n-2}}{K_{m-1}} by \lemref{TwoCliquesCharac}. Thus
$\beta(\CP{K_{n-2}}{K_{m-1}})\leq|S'|\leq|S|-2=\beta(\CP{K_n}{K_m})-2$, which implies \eqref{eqn:TwoCliquesProductB}. }


\BLAH{Proof of \thmref{TwoCliquesProduct}.}{We proceed by induction on $n+m$ in increments of $3$. (Formally speaking, we are doing induction on $\floor{\tfrac{1}{3}(n+m)}$.)\

First observe that for $m=1$, we know that
$\beta(\CP{K_n}{K_m})=n-1$. For $m=2$, we have
$\beta(\CP{K_2}{K_2})=2=\floor{\twothirds(2+2-1)}$,
$\beta(\CP{K_3}{K_2})=2=\floor{\twothirds(3+2-1)}$, and
$\beta(\CP{K_n}{K_2})=n-1$ for all $n\ge 3$. Thus the assertion is
true for $m\leq2$. Now suppose that $m\geq3$. By
\lemref{TwoCliquesProduct} we have
\EQN{FromLemma}{\beta(\CP{K_n}{K_m})=
2+\min\{\beta(\CP{K_{n-2}}{K_{m-1}}),\beta(\CP{K_{n-1}}{K_{m-2}})\}.}

\textbf{Case 1.} $n\geq2m-1$: Then $n\geq2\cdot\beta(K_m)+1$ by
\eqnref{BetaCompleteGraph}, and $\beta(\CP{K_n}{K_m})=n-1$ by
\thmref{BigClique} with $G=K_m$.

\textbf{Case 2.} $n=2m-2$: First consider \CP{K_{n'}}{K_{m'}},
where $n'=n-1=2m-3$ and $m'=m-2$. Then $m'\le n'$ and $n'\ge
2m'-1$. By induction, \begin{equation*}
\beta(\CP{K_{n'}}{K_{m'}})=n'-1=n-2=\floor{\twothirds(n+m-1)}-2.
\end{equation*}
Now consider \CP{K_{n'}}{K_{m'}}, where $m'=m-1$ and
$n'=n-2=2m-4$. Then $m'\le n'\le 2m'-1$. By induction
\begin{equation*}
\beta(\CP{K_{n'}}{K_{m'}})=\floor{\twothirds(n'+m'-1)}
=\floor{\twothirds(n+m-1)}-2.
\end{equation*}
By \eqnref{FromLemma}, $\beta(\CP{K_n}{K_m})=\floor{\twothirds(n+m-1)}$.

\textbf{Case 3.} $n=2m-3$: First consider \CP{K_{n'}}{K_{m'}},
where $m'=m-2$ and $n'=n-1=2m-4$. Then $m'\le n'$ and $n'\ge
2m'-1$. By induction, \begin{equation*}
\beta(\CP{K_{n'}}{K_{m'}})=n'-1=n-2=\floor{\twothirds(n+m-1)}-2.
\end{equation*}
Now consider \CP{K_{n'}}{K_{m'}}, where $m'=m-1$, $n'=n-2=2m-5$.
For $m\ge 4$, we have $m'\le n' \le 2m'-1$. By induction
\begin{equation*} \beta(\CP{K_{n'}}{K_{m'}})=
\floor{\twothirds(n'+m'-1)}=\floor{\twothirds(n+m-1)}-2.
\end{equation*}
For $m=3$, we have $n=2m-3=3$. It is easily verified that
$\beta(\CP{K_3}{K_3})=3=\floor{\twothirds(3+3-1)}$. In all cases
we obtain $\beta(\CP{K_n}{K_m})=\floor{\twothirds(n+m-1)}$ by
\eqnref{FromLemma}.

\textbf{Case 4.} $n\le 2m-4$:  First consider \CP{K_{n'}}{K_{m'}}, where $m'=m-2$ and $n'=n-1\le 2m-5$. Then,
$m'\le n'\le 2m'-1$. By induction, \begin{equation*}
\beta(\CP{K_{n'}}{K_{m'}})=\floor{\twothirds(n'+m'-1)}= \floor{\twothirds(n+m-1)}-2.
\end{equation*}
Now consider \CP{K_{n'}}{K_{m'}}, where $m'=m-1$ and $n'=n-2\le 2m-6$.
If $m\le n-1$, then $m'\le n'< 2m'-1$, and by induction
\begin{equation*}
\beta(\CP{K_{n'}}{K_{m'}})=\floor{\twothirds(n'+m'-1)}=
\floor{\twothirds(n+m-1)}-2.
\end{equation*}
If $m=n\ge 4 $,  then $n'\le m' \le 2n'-1$ and by
induction
\begin{equation*}
\beta(\CP{K_{m'}}{K_{n'}})=\floor{\twothirds(m'+n'-1)}= \floor{\twothirds(n+m-1)}-2.
\end{equation*}
Finally, if $m=n=3$,  then $\beta(\CP{K_{n'}}{K_{m'}})=\beta(\CP{K_2}{K_1})=1=\floor{\twothirds(3+3-1)}
-2$. In all cases, we obtain $\beta(\CP{K_n}{K_m})=\floor{\twothirds(m+n-1)}$ by \eqnref{FromLemma}.}

\SECT{Paths}{Paths and Grids}

Let $P_n$ denote the path on $n\geq1$ vertices. \citet{KRR-DAM96}
and \citet{CEJO-DAM00} proved that an $n$-vertex graph $G$ has
\EQN{BetaPath}{\beta(G)=1\,\Longleftrightarrow\,G=P_n.} Thus, by
\thmref{BigClique}, for all $n\geq3$,
\EQN{CliquePath}{\beta(\CP{K_n}{P_m})=n-1.} Minimum doubly
resolving sets in paths are easily characterised.

\LEM{DoublyResolvePath}{For all $n\geq2$ we have $\psi(P_n)=2$.
Moreover, the two endpoints of $P_n$ are in every doubly resolving
set of $P_n$.}

\PROOF{By definition $\psi(G)\geq2$ for every graph $G\ne K_1$.
Let $P_n=(v_1,v_2,\dots,v_n)$. For all $1\leq i<j\leq n$, we have
$d(v_i,v_1)-d(v_j,v_1)=(i-1)-(j-1)=i-j$, and
$d(v_i,v_n)-d(v_j,v_n)=(n-i)-(n-j)=j-i$. Thus $\{v_1,v_n\}$ doubly
resolve $P_n$, and $\psi(P_n)=2$. Finally, observe that $v_1$ is
in every doubly resolving set, as otherwise $v_1$ and $v_2$ would
not be doubly resolved. Similarly $v_n$ is in every doubly
resolving set.}

\LEM{BetaProductTwo}{If $\beta(\CP{G}{H})=2$, then $G$ or $H$ is a path.}

\PROOF{Say $S=\{av,bw\}$ resolves \CP{G}{H}. Suppose that $a=b$.
Then the projection of $S$ onto $G$ is a single vertex. By
\lemref{Projections}, the projection of $S$ onto $G$ resolves $G$,
and by \eqnref{BetaPath}, only paths have singleton resolving
sets. Thus $G$ is a path, and we are done. Similarly, if $v=w$
then $H$ is a path, and we are done.  Now suppose that $a\ne b$
and $v\ne w$. Let $c$ be the neighbour of $b$ on a shortest path
from $a$ to $b$. Note that $c$ may equal $a$. Then
$d_G(a,c)+1=d_G(a,b)$ and $d_G(b,c)=1$. Similarly, let $x$ be the
neighbour of $w$ on a shortest path from $v$ to $w$. Then
$d_H(v,x)+1=d_H(v,w)$ and $d_H(x,w)=1$. This implies that $S$ does
not resolve $bx$ and $cw$, since
\begin{equation*}
d(bx,av)=d_G(a,b)+d_H(x,v)=d_G(a,c)+d_H(v,w)=d(cw,av)
\end{equation*}
and
\begin{equation*}
d(bx,bw)=d_H(x,w)=1=d_G(b,c)=d(cw,bw).
\end{equation*}
This contradiction proves the result.}

\thmref{DoublyResolving} and \lemref{DoublyResolvePath} imply that
every graph $G$ satisfies
\EQN{MultiplyPath}{\beta(G)\leq\beta(\CP{G}{P_n})\leq\beta(G)+1,}
as proved by \citet{CEJO-DAM00} in the case that $n=2$.

An $n$-dimensional \emph{grid} is a cartesian product of paths
$\CCP{P_{m_1}}{P_{m_2}}{P_{m_n}}$.
\twoeqnref{BetaPath}{MultiplyPath} imply that,
\EQN{Grid}{\beta(\CCP{P_{m_1}}{P_{m_2}}{P_{m_n}})\leq n.} as
proved by \citet{KRR-DAM96}, who in addition claimed that
$$\beta(\CCP{P_{m_1}}{P_{m_2}}{P_{m_n}})=n.$$ They wrote `we leave
it for the reader to see why $n$ is a lower bound'. This claim is
false if every $m_i=2$ and $n$ is large, since
$\beta(\CCP{P_2}{P_2}{P_2})\rightarrow 2n/\log n$ as discussed in
\secref{Coin}. \citet{ST-MOR04} claimed without proof that `using
a result of \citet{Lindstrom-CMB65}' one can prove that
\EQN{Claim}{\limsup_{n\rightarrow\infty}\,\beta(\underbrace{\CCP{P_k}{P_k}{P_k}}_n)\cdot\frac{\log
n}{n\log k}\leq2.}

\SECT{Cycles}{Cycles}

Let $C_n$ denote the cycle on $n\geq3$ vertices. Two vertices $v$
and $w$ of $C_n$ are \emph{antipodal} if $d(v,w)=\frac{n}{2}$.
Note that no two vertices are antipodal in an odd cycle.

\LEMMA{\citep{KRR-DAM96,SZ-CMJ04}}{BetaCycle}{For all $n\geq3$ we
have $\beta(C_n)=2$. Moreover, two vertices resolve $C_n$ if and
only if they are not antipodal.}

\LEM{PsiCycle}{For all $n\geq3$ we have
\begin{equation*}
\psi(C_n)=
\begin{cases}
2   & \textup{, if }n\textup{ is odd}\\
3   & \textup{, if $n$ is even.}\\
\end{cases}
\end{equation*}}

\PROOF{We have $\psi(C_n)\geq2$ by definition. Now we prove the
upper bound. Denote $C_n=(v_1,v_2,\dots,v_n)$. Let
$k:=\floor{\frac{n}{2}}$. Consider two vertices $v_i$ and $v_j$ of
$C_n$. Without loss of generality $i<j$.

\textbf{Case 1.} $1\leq i<j\leq k+1$: Then
$d(v_i,v_1)-d(v_j,v_1)=(i-1)-(j-1)=i-j$, and
$d(v_i,v_{k+1})-d(v_j,v_{k+1})=(k+1-i)-(k+1-j)=j-i\ne i-j$.
Thus $v_1,v_{k+1}$ doubly resolve $v_i,v_j$.

\textbf{Case 2.} $k+1\leq i<j\leq n$: Then
$d(v_i,v_1)-d(v_j,v_1)=(n+1-i)-(n+1-j)=j-i$, and
$d(v_i,v_{k+1})-d(v_j,v_{k+1})=(i-k-1)-(j-k-1)=i-j\ne j-i$.
Thus $v_1,v_{k+1}$ doubly resolve $v_i,v_j$.

\textbf{Case 3.} $1\leq i\leq k+1<j\leq n$: Suppose that
$v_1,v_{k+1}$ does not doubly resolve $v_i,v_j$.
That is, $d(v_i,v_1)-d(v_j,v_1)=d(v_i,v_{k+1})-d(v_j,v_{k+1})$.
Thus $(i-1)-(n+1-j)=(k+1-i)-(j-k-1)$.
Hence $n=2i+2j-2k-4$ is even.

Therefore for odd $n$, $\{v_1,v_{k+1}\}$ doubly resolves $C_n$, and $\psi(C_n)=2$.

For even $n$, in Case 3, suppose that $v_1,v_2$ does not doubly
resolve $v_i,v_j$. That is,
$d(v_i,v_1)-d(v_j,v_1)=d(v_i,v_2)-d(v_j,v_2)$. Thus
$(i-1)-(n+1-j)=(i-2)-(n+2-j)$ and $-2=-4$, a contradiction. Hence
for even $n$, $\{v_1,v_2,v_{k+1}\}$ doubly resolve $C_n$, and
$\psi(C_n)\leq3$.

It remains to prove that $\psi(C_n)\geq3$ for even $n$. Suppose
that $\psi(C_n)\leq2$ for some even $n=2k$. By symmetry we can
assume that $\{v_1,v_i\}$ doubly resolves $C_n$ for some $2\leq
i\leq k+1$.

\textbf{Case 1.} $2\leq i\leq k-1$: Then $d(v_{i+1},v_1)-d(v_{i+2},v_1)=i-(i+1)=-1$, and
$d(v_{i+1},v_i)-d(v_{i+2},v_i)=1-2=-1$. Thus $v_1,v_i$ does not resolve $v_{i+1},v_{i+2}$.

\textbf{Case 2.} $i=k$: Then $d(v_2,v_1)-d(v_{n-1},v_1)=1-2=-1$, and
$d(v_2,v_i)-d(v_{n-1},v_i)=(k-2)-(k-1)=-1$. Thus $v_1,v_i$ does not resolve $v_2,v_{n-1}$.

\textbf{Case 3.} $i=k+1$: Then $d(v_2,v_1)-d(v_n,v_1)=1-1=0$, and
$d(v_2,v_i)-d(v_n,v_i)=(k-1)-(k-1)=0$. Thus $v_1,v_i$ does not resolve $v_2,v_n$.

In each case we have derived a contradiction. Thus $\psi(C_n)\geq3$ for even $n$.}


\thmref{DoublyResolving} and \lemref{PsiCycle} imply that
every graph $G$ satisfies
\EQN{GraphCycle}{\beta(G)\leq\beta(\CP{G}{C_n})\leq
\begin{cases}
\beta(G)+1  & \textup{, if }n\textup{ is odd}\\
\beta(G)+2  & \textup{, if }n\textup{ is even.}
\end{cases}}


\THM{GraphCycle}{For every graph $G$ and for all $n\geq3$, we have
$\beta(\CP{G}{C_n})=2$ if and only if $G$ is a path and $n$ is
odd.}

\PROOF{($\Longleftarrow$) Since $G$ is a path, $\beta(G)=1$ by
\eqnref{BetaPath}. Since $n$ is odd, $\psi(C_n)=2$ by
\lemref{PsiCycle}. Thus
$\beta(\CP{G}{C_n})\leq\psi(C_n)+\beta(G)-1=2$ by
\thmref{DoublyResolving}.

($\Longrightarrow$) Suppose that $\beta(\CP{G}{C_n})=2$. Say
$S=\{av,bw\}$ resolves \CP{G}{C_n}. Then $G$ is a path by
\lemref{BetaProductTwo}. It remains to show that $n$ is odd.
Suppose on the contrary that $n=2r$ is even. Let $C=C_n$. By
\corref{Projections}, the projection $\{v,w\}$ of $S$ onto $C$
resolves $C$. By \lemref{BetaCycle}, we have $\beta(C)=2$, and
thus $v\ne w$. Moreover, $v$ and $w$ are not antipodal. That is,
$d_C(v,w)\leq r-1$. Hence there is a neighbour $x$ of $w$ in $C$
with $d_C(v,x)=d_C(v,w)+1$. Now consider $G$. If $a\ne b$, then
using the argument from the proof of \lemref{BetaProductTwo}, we
can construct a pair of vertices that are not resolved by $S$. So
now assume $a=b$. That is, our resolving set is contained in a
single column of \CP{G}{C_n}. Let $p$ be a neighbour of $a$ in
$G$. Then $S$ does not resolve $pw$ and $ax$, since
$d(pw,bw)=1=d(ax,bw)$ and $d(pw,av)=1+d_C(v,w)=d_C(x,v)=d(ax,av)$.
This contradiction proves the result.}


By \lemref{PsiCycle} and \eqnref{BetaPath}, we have
$\beta(\CP{P_m}{C_n})\leq\psi(C_n)+\beta(P_m)-1\leq 3+1-1=3$. Thus
\thmref{GraphCycle} implies that for all $m\geq2$ and $n\geq3$ we
have \EQN{PathCycle}{ \beta(\CP{P_m}{C_n})=
\begin{cases}
2   & \textup{, if }n\textup{ is odd}\\
3   & \textup{, if }n\textup{ is even.}
\end{cases}}

\THM{CycleCycle}{For all $m,n\geq3$ we have
$$\beta(\CP{C_m}{C_n})=
\begin{cases}
3   & \textup{, if }m\textup{ or }n\textup{ is odd}\\
4   & \textup{, otherwise}.
\end{cases}$$
}

\PROOF{We have $\beta(\CP{C_m}{C_n})\geq3$ by \thmref{GraphCycle}.
If $m$ or $n$ is odd, then $\beta(\CP{C_m}{C_n})\leq3$ by
\eqnref{GraphCycle} and since $\beta(C_m)=2$. It remains to prove
that $\beta(\CP{C_m}{C_n})\geq4$ when $m$ and $n$ are even. Let
$G:=\CP{C_{2r}}{C_{2s}}$. We denote each vertex $U$ of $G$ by
$u_1u_2$, where $u_1\in C_{2r}$ and $u_2\in C_{2s}$.

Observe that in $C_{2r}$, every vertex $u$ is antipodal with a
unique vertex $v$; thus $d(x,u)+d(x,v)=r$ for every vertex $x$ of
$C_{2r}$.

Two vertices $U$ and $V$ of $G$ are \emph{antipodal} if $u_1$ and
$v_1$ are antipodal in $C_{2r}$ and $u_2$ and $v_2$ are antipodal
in $C_{2s}$. Suppose that $U$ and $V$ are antipodal. Then for
every vertex $W$ of $G$, we have
\EQN{jhjhjh}{d_G(W,U)+d_G(W,V)=d(w_1,u_1)+d(w_2,u_2)+d(w_1,v_1)+d(w_2,v_2)=r+s.}

\begin{claim}
Let $U$ be a vertex in a resolving set $S$ of $G$. Say $U$ and $V$
are antipodal. Then the set $S'$ obtained by replacing $U$ by $V$
in $S$ also resolves $G$.
\end{claim}

\PROOF{Suppose on the contrary, that $S'$ does not resolve $G$.
Thus there exist vertices $X,Y$ of $G$ such that
$d_G(X,Z)=d_G(Y,Z)$ for every vertex $Z\in S'$. In particular,
$d_G(X,V)=d_G(Y,V)$. By \eqnref{jhjhjh},
$d_G(X,U)-r-s=d_G(Y,U)-r-s$, implying $d_G(X,U)=d_G(Y,U)$. Thus
$d_G(X,Z)=d_G(Y,Z)$ for every vertex $Z\in S$; that is, $X$ and
$Y$ are not resolved by $S$. This contradiction proves the claim.}

Suppose on the contrary that $S=\{ U,V,W \}$ is a resolving set of
$G$. Represent $G$ by the points of a $2r\times 2s$ grid.
Consecutive points in the same row or column are adjacent, and the
first and last points of the same row or column are adjacent.
Observe that antipodal vertices of $G$ are in opposite quadrants
of the grid. Thus, by the above claim, we can assume that $U,V,W$
are in one of the four halves of the grid. Without loss of
generality, $U,V,W$ are in the left half of the grid. This implies
that $d(u_1,v_1)<r$, $d(u_1,w_1)<r$ and $d(v_1,w_1)<r$.
Furthermore, $U,V,W$ are in at least two different rows and two
different columns, since the projections of $S$ resolve $C_{2r}$
and $C_{2s}$.

By symmetry, it suffices to consider the following cases:\\
1. $U,V,W$ are in different rows and different columns,\\
2. $U,V,W$ are in different rows, but $U,V$ are in the same column, and\\
3. $U,V$ are in the same column and $V,W$ in the same row. \\
In each case we will find vertices $X,Y$ such that
$d(X,U)=d(Y,U)$, $d(X,V)=d(Y,V)$ and $d(X,W)=d(Y,W)$; that is, $S$
does not resolve the pair $X,Y$.

\textbf{Case 1.} Assume that, if one of the vertices $u_2,v_2,w_2$
is in the shortest path determined by the other two vertices, then
that vertex is $v_2$. It is then possible to draw the grid in such
a way that the projections $u_2,v_2,w_2$ appear from bottom to top
in $C_{2s}$, $d(u_2,v_2)<s$, and $d(v_2,w_2)<s$. Now, if $v_1$ is
in the shortest path between $u_1$ and $w_1$ in $C_{2r}$, then let
$X,Y$ be the two neighbours of $V$ lying in shortest paths between
$V$ and $W$; see \figref{I}(a). Otherwise, assume that $u_1$ is in
the shortest path between $v_1$ and $w_1$.  Let $Z$ be the vertex
$u_1v_2$. Let $X,Y$ be the neighbours of $Z$ in shortest paths
between $Z$ and $W$; see \figref{I}(b). It is easy to verify that
in both cases $d(X,U)=d(Y,U)$, $d(X,V)=d(Y,V)$ and
$d(X,W)=d(Y,W)$.

\begin{figure}[!ht]
\begin{center}
\includegraphics{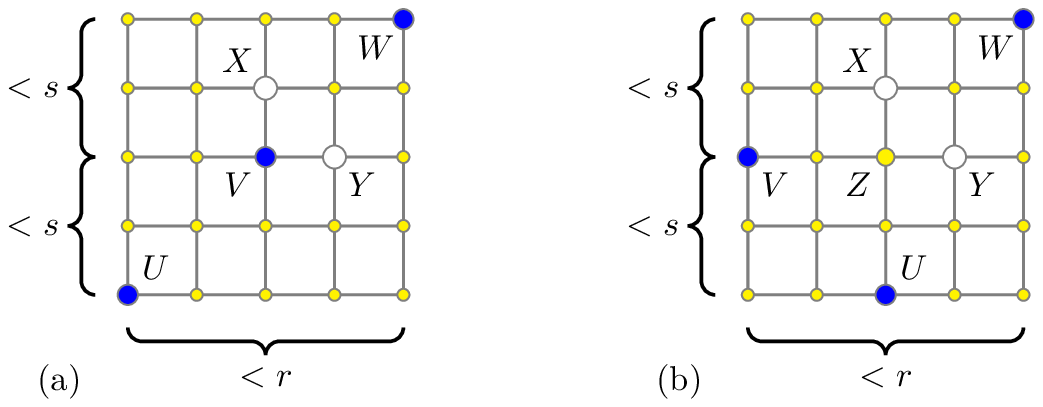}
\caption{\label{fig:I}Illustration for Case 1 of \thmref{CycleCycle}.}
\end{center}
\end{figure}

\textbf{Case 2.} Observe that at least two of the distances
$d(u_2,v_2)$, $d(v_2,w_2)$ and $d(u_2,w_2)$ in $C_{2s}$ must be
less than $s$. If $u_2,v_2$ are not antipodal in $C_{2s}$ and
$w_2$ is not in the shortest path between $u_2$ and $v_2$ in
$C_{2s}$, then $d(u_2,w_2)<s$ or $d(v_2,w_2)<s$. Let us assume
that $d(v_2,w_2)<s$. Let $X,Y$ be the vertices adjacent to $V$
lying in a shortest path between $V$ and $W$; see \figref{II}(a).
If $u_2,v_2$ are not antipodal in $C_{2s}$ and $w_2$ is in the
shortest path between $u_2$ and $v_2$ in $C_{2s}$, then let $X,Y$
be the neighbours of $V$ not lying in a shortest path between $V$
and $W$; see \figref{II}(b). Finally, if $u_2,v_2$ are antipodal
in $C_{2s}$, consider the vertices $X,Y$ at distance two from $V$;
see \figref{II}(c). It is easy to verify that in all cases
$d(X,U)=d(Y,U)$, $d(X,V)=d(Y,V)$ and $d(X,W)=d(Y,W)$.

\begin{figure}[!ht]
\begin{center}
\includegraphics{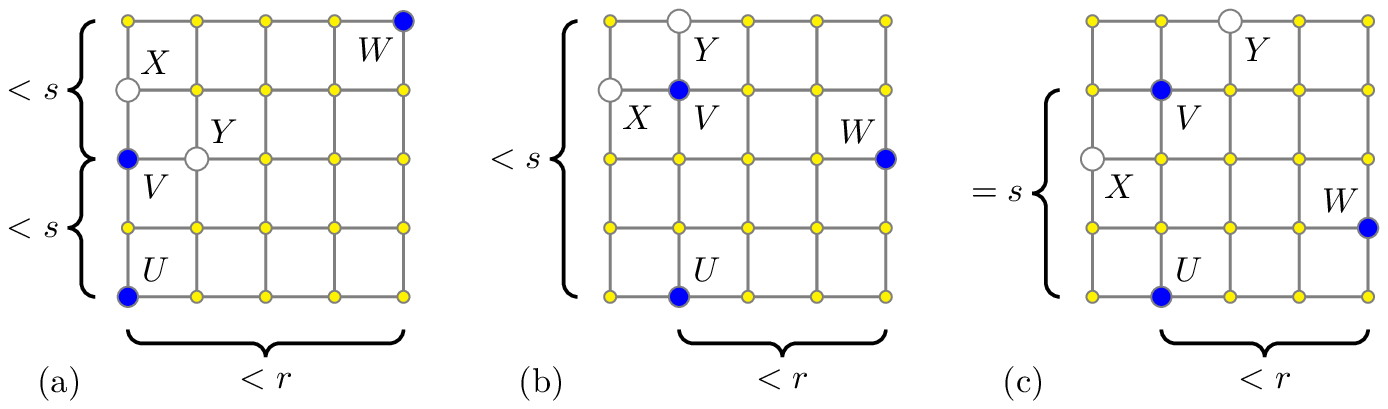}
\caption{\label{fig:II}Illustration for Case 2 of \thmref{CycleCycle}.}
\end{center}
\end{figure}

\begin{figure}[!ht]
\begin{center}
\includegraphics{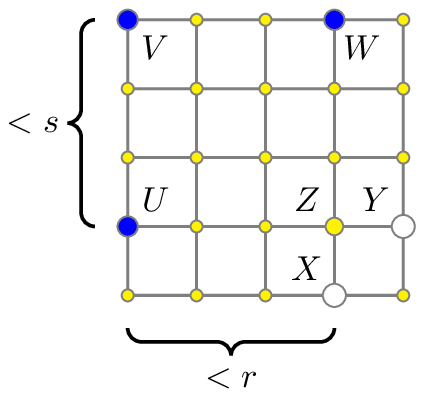}
\caption{\label{fig:III}Illustration for Case 3  of \thmref{CycleCycle}.}
\end{center}
\end{figure}

\textbf{Case 3.} In this case, $d(u_2,v_2)<s$ since the projection
$\{u_2,v_2,w_2\}=\{u_2,v_2\}$ resolves $C_{2s}$. Let
$Z:=(w_1,u_2)$. Let $X,Y$ be the neighbours of $Z$ not lying in a
shortest path between $Z$ and $V$; see \figref{III}. It is easy to
verify that $d(X,U)=d(Y,U)$, $d(X,V)=d(Y,V)$ and $d(X,W)=d(Y,W)$.}


\THM{CompleteCycle}{For all $n\geq1$ and $m\geq3$ we have
\begin{equation*}
\beta(\CP{K_n}{C_m})=
\begin{cases}
2   & \textup{, if $n=1$}\\
2   & \textup{, if $n=2$ and $m$ is odd}\\
3   & \textup{, if $n=2$ and $m$ is even}\\
3   & \textup{, if $n=3$}\\
3   & \textup{, if $n=4$ and $m$ is even}\\
4   & \textup{, if $n=4$ and $m$ is odd}\\
n-1 & \textup{, if $n\geq5$.}\\
\end{cases}
\end{equation*}}

\PROOF{The case $n\geq2\beta(C_n)+1=5$  is an immediate corollary
of \thmref{BigClique} and \lemref{BetaCycle}. The case $n=3$ is a
special case of \thmref{CycleCycle} since $K_3=C_3$. The case
$n=2$ is a special case of \eqnref{PathCycle} since $K_2=P_2$. The
case $n=1$ is a repetition of \lemref{BetaCycle}.

It remains to prove the case $n=4$. Say $V(K_4)=\{a,b,c,d\}$.
First note that $\beta(\CP{K_4}{C_m})\geq\beta(K_4)=3$ by
\corref{Projections} and \eqnref{BetaCompleteGraph}. By
\lemref{PsiCompleteGraph} we have $\psi(K_4)=3$. Thus
$\beta(\CP{K_4}{C_m})\leq4$ by \lemref{BetaCycle} and
\thmref{DoublyResolving} with $H=K_4$. For even $m$, it is easily
verified that $\{av,bv,cw\}$ resolves \CP{K_4}{C_m} for any edge
$vw$ of $C_m$.


It remains to prove that $\beta(\CP{K_4}{C_m})\geq4$ for odd
$m=2h+1$. Consider the vertices of \CP{K_4}{C_m} to be in a
$4\times m$ grid, where two vertices in the same row are adjacent,
and two vertices in the same column are adjacent if and only if
they are consecutive rows or they are in the first and last rows.
Suppose on the contrary that $S=\{u,v,w\}$ resolves \CP{K_4}{C_m}.
Then $u,v,w$ are in three different columns and in at least two
different rows (by considering the projections of $S$ onto $K_4$
and $C_m$).

\textbf{Case 1.} Suppose that two vertices in $S$, say $u$ and
$v$, are in the same row. Consider the grid centred at the row of
$u,v$. Without loss of generality, $u$ and $v$ are in the first
and second columns, and $w$ is in a row above $u$ and $v$. Let $x$
and $y$ be the vertices shown in \figref{CompleteCycle}(a). Then
$d(x,u)=d(y,u)=h+1$, $d(x,v)=d(y,v)=h+1$, and $d(x,w)=d(y,w)=p$.
Thus $S$ does not resolve $x$ and $y$, which is the desired
contradiction.

\begin{figure}[!ht]
\begin{center}
\includegraphics{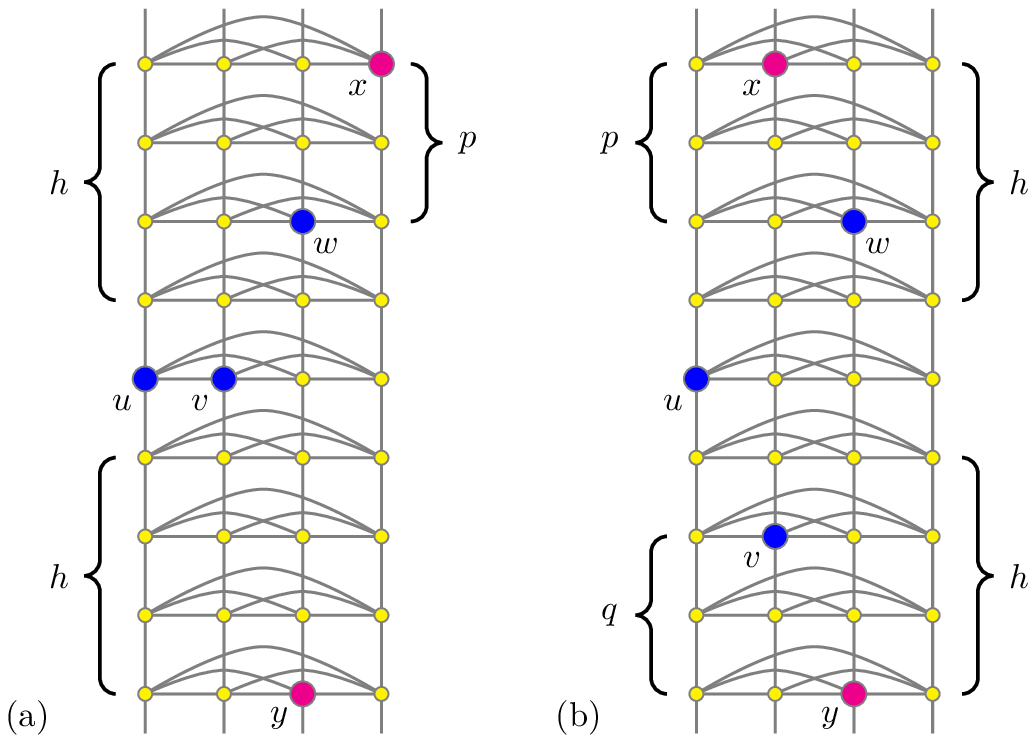}
\caption{\label{fig:CompleteCycle}Illustration for \thmref{CompleteCycle}.}
\end{center}
\end{figure}

\textbf{Case 2.} Now suppose that $u,v,w$ are in different rows.
Without loss of generality, $u$ is in the middle row and the first
column, and  $v$ is in the second column and in a row below $u$,
and $w$ is in the third column and in a row above $u$. Let $x$ and
$y$ be the vertices shown in \figref{CompleteCycle}(b). Then
$d(x,u)=d(y,u)=h+1$, $d(x,v)=d(y,v)=q$, and $d(x,w)=d(y,w)=p$.
Thus $S$ does not resolve $x$ and $y$, which is the desired
contradiction.}

%

\SECT{Trees}{Trees}

Let $v$ be a vertex of a tree $T$. Let $\ell_v$ be the number of components of $T\setminus v$ that are (possibly edgeless) paths. \citet{KRR-DAM96} and \citet{CEJO-DAM00} proved that for every tree $T$ that is not a path,
\EQN{BetaTree}{\beta(T)=\sum_{v\in V(T)}\max\{\ell_v-1,0\}.}
A \emph{leaf} of a graph is a vertex of degree one. The following result
for doubly resolving sets in trees is a generalisation of \lemref{DoublyResolvePath} for paths.

\LEM{DoublyResolvingTree}{The set of leaves $L$ is the unique minimum doubly resolving set for a tree $T$, and $\psi(T)=|L|$.}

\PROOF{Every pair of vertices $v,w$ of $T$ lie on a path whose endpoints
are leaves $x,y$. Clearly $x,y$ doubly resolve $v,w$. Thus $L$ is a
doubly resolving set. Say $v$ is a leaf of $T$ whose neighbour is
$w$. Every shortest path from $v$ passes through $w$.
Thus $v,w$ can only be doubly resolved by a pair including $v$.
Thus $v$ is in every doubly resolving set of $T$. The result follows.}

\thmref{DoublyResolving} and \lemref{DoublyResolvingTree} imply that for every tree $T$ with $k$ leaves and for every graph $G$,
\begin{equation}
\beta(\CP{T}{G})\leq\beta(G)+k-1.
\end{equation}
Moreover, many leaves force up the metric dimension of a cartesian product.

\LEM{LeavesLowerBound}{Every graph $G$ with $k\geq2$ leaves satisfies $\beta(\CP{G}{G})\geq k$.}

\PROOF{Let $S$ with a metric basis of \CP{G}{G}. Let $b$ and $w$ be distinct leaves of $G$ respectively adjacent to $a$ and $v$. There is a vertex $xy\in S$ that resolves $aw$ and $bv$. Suppose on the contrary that $x\ne b$ and $y\ne w$.
Thus $d_G(b,x)=d_G(a,x)+1$ and $d_G(w,y)=d_G(v,y)+1$. Hence $d_G(a,x)-d_G(b,x)=d_G(v,y)-d_G(w,y)=-1$, which implies that $d_G(a,x)+d_G(w,y)=d_G(b,x)+d_G(v,y)$. That is, $d(aw,xy)=d(bv,xy)$.
Thus $xy$ does not resolve $aw$ and $bv$. This contradiction proves that $x=b$ or $y=w$. Thus for every pair of leaves $b,w$ there is a vertex $by$ or $xw$ in
$S$. Suppose that for some leaf $b$, there is no vertex $by\in S$. Then for every leaf $w$, there is a vertex $xw\in S$, and $|S|\geq k$. Otherwise for every leaf $b$, there is a vertex $by\in S$, and again $|S|\geq k$.}

The following result implies that $\psi$ is not bounded by any function of metric dimension.

\THM{Comb}{For every integer $n\geq4$ there is a tree $B_n$ with $\beta(B_n)=2$ and
\begin{equation*}
n=\psi(B_n)\leq\beta(\CP{B_n}{B_n})\leq n+1.
\end{equation*}}

\begin{figure}[!ht]
\begin{center}
\includegraphics{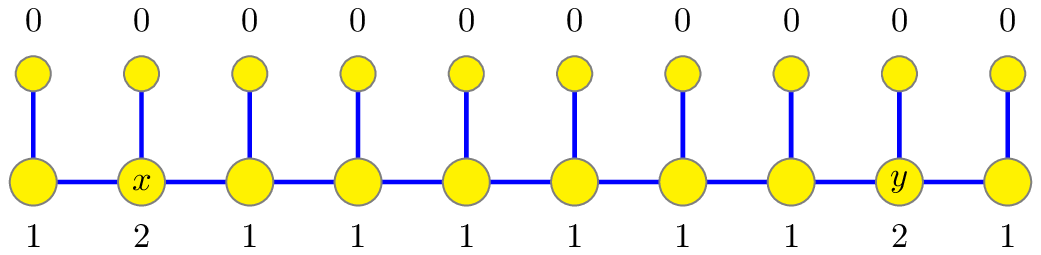}
\caption{\label{fig:Comb}An illustration  of the comb graph $B_{10}$ showing the $\ell$-values at each vertex.}
\end{center}
\end{figure}

\PROOF{Let $B_n$ be the \emph{comb} graph obtained by attaching one leaf at
every vertex of $P_n$. Now $\ell_v=0$ for every leaf $v$ of $B_n$, and
$\ell_w=1$ for every other vertex $w$ of $B_n$, except for the two vertices $x$ and $y$ indicated in \figref{Comb}, for which $\ell_x=\ell_y=2$. Thus $\beta(B_n)=2$ by \eqnref{BetaTree}. Since $B_n$ has $n$ leaves, we have $\psi(B_n)=n$ by \lemref{DoublyResolvingTree}. Moreover, $\beta(\CP{B_n}{B_n})\geq n$ by \lemref{LeavesLowerBound}. The upper bound $\beta(\CP{B_n}{B_n})\leq n+1$ follows from \thmref{DoublyResolving}.}

Given that the proof of \thmref{Comb} is heavily dependent on the presence of leaves in $B_n$, it is tempting to suspect that such behaviour does not occur among more highly connected graphs. This is not the case.

\THM{HighlyConnected}{For all $k\geq1$ and $n\geq2$ there is a $k$-connected graph $G_{n,k}$ for which $\beta(G_{n,k})\leq 2k$ and $\beta(\CP{G_{n,k}}{G_{n,k}})\geq n$.}

\begin{figure}[!ht]
\begin{center}
\includegraphics{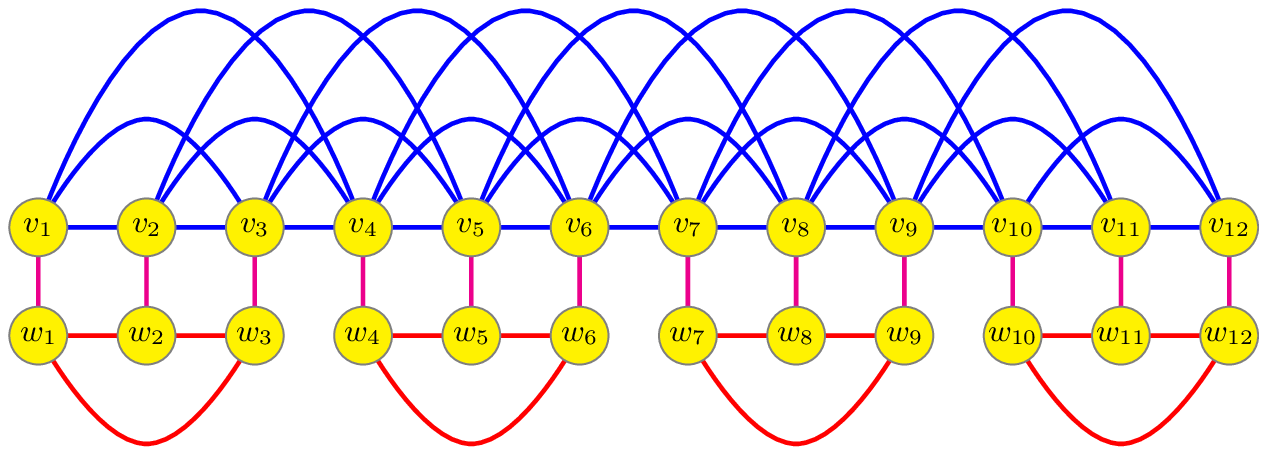}
\caption{\label{fig:Example}The construction in \thmref{HighlyConnected} with $k=3$ and $n=2$.}
\end{center}
\end{figure}

\PROOF{As illustrated in \figref{Example}, let $G_{n,k}$ be the graph with vertex set $\{v_i,w_i:1\leq i\leq 2kn\}$, where every $v_iw_i$ is an edge, $v_iv_j$ is an edge whenever $|i-j|\leq k$, and $w_iw_j$ is an edge whenever $\ceil{i/k}=\ceil{j/k}$. Note that $G_{n,1}=B_{2n}$. Clearly $G_{n,k}$ is $k$-connected. It is easily seen that $\{v_i,v_{2kn+1-i}:1\leq i\leq k\}$ resolves $G_{n,k}$. Thus $\beta(G_{n,k})\leq 2k$.

Say $S$ doubly resolves $G_{n,k}$. On the contrary, suppose that $$S\cap\{w_{\ell k+1},w_{\ell k+2},\dots,w_{\ell k+k}\}=\emptyset$$ for some $\ell$ with $0\leq\ell\leq2n-1$. This implies that $d(w_{\ell k+1},x)=d(v_{\ell k+1},x)+1$ for every vertex $x\in S$. Hence $S$ does not doubly resolve $w_{\ell k+1}$ and $v_{\ell k+1}$. This contradiction proves that $S\cap\{w_{\ell k+1},w_{\ell k+2},\dots,w_{\ell k+k}\}\ne\emptyset$ for every $\ell $ with $0\leq\ell\leq2n-1$. Thus $|S|\geq2n$ and $\psi(G_{n,k})\geq2n$. That $\beta(\CP{G_{n,k}}{G_{n,k}})\geq n$ follows from \lemref{TwoProjections}.}

We conclude that for all $k\geq1$, there is no function $f$ such that
$\beta(\CP{G}{H})\leq f(\beta(G),\beta(H))$ for all $k$-connected graphs $G$ and $H$.

\def\soft#1{\leavevmode\setbox0=\hbox{h}\dimen7=\ht0\advance \dimen7
  by-1ex\relax\if t#1\relax\rlap{\raise.6\dimen7
  \hbox{\kern.3ex\char'47}}#1\relax\else\if T#1\relax
  \rlap{\raise.5\dimen7\hbox{\kern1.3ex\char'47}}#1\relax \else\if
  d#1\relax\rlap{\raise.5\dimen7\hbox{\kern.9ex \char'47}}#1\relax\else\if
  D#1\relax\rlap{\raise.5\dimen7 \hbox{\kern1.4ex\char'47}}#1\relax\else\if
  l#1\relax \rlap{\raise.5\dimen7\hbox{\kern.4ex\char'47}}#1\relax \else\if
  L#1\relax\rlap{\raise.5\dimen7\hbox{\kern.7ex
  \char'47}}#1\relax\else\message{accent \string\soft \space #1 not
  defined!}#1\relax\fi\fi\fi\fi\fi\fi} \def\cprime{$'$}

\end{document}